\documentclass{amsart}

\theoremstyle{definition}

\theoremstyle{remark}

\numberwithin{equation}{section}

\begin{document}

\title{Trading Strategy Adopted Optimization of European Call Option}

\author{Toshio Fukumi}
\address{Department of Business Administration, Matsuyama University, 2-4 Bunkyo, Matsuyama, Ehime 790 8578, Japan}
\email{tofukumi@cc.matsuyama-u.ac.jp}



\date{}



\begin{abstract}
Optimal pricing of European call option is described by linear stochastic differential equation with mean income and volatility. Trading strategy is given by a twin of stochastic variables for which self-financial condition is not postulated. Optimal pricing will be deformed to be adoptive to trading strategy employing martingale property where stochastic integral w.r.t. analytical solution of Black-Scholes type stochastic differential equation  will be employed.
\end{abstract}

\maketitle


\par

\section {Stochastic Optimization}
.
\subsection{Black-Scholes Formula}
\par
Black-Scholes' formula starts from heat equation
\begin{equation}
{\partial u\over {\partial t}}={1\over 2} {\partial^2 u^2\over {\partial^2 x^2}},
\end{equation}
the solution is given by
\begin{equation}
 u= {1 \over {\sqrt {2{\pi}t}}}e^  {-{u^2 \over 2t}},
\end{equation}
and corresponding stochastic differential equation is
\begin{equation}
 du_t = dB_t,
\end{equation}
where $B_t$ denotes Brownian motion. In case there is a drift term we have
\begin{equation}
 du_t =  -cdt+dB_t,
\end{equation}
whose solution is
\begin {equation}
 u= {1 \over {\sqrt {2{\pi}t}}}e^  {-{(u-ct)^2 \over 2t}},
\end{equation}
where $c$ is a velocity.
\par
Let us modify the last equation to be
\begin {equation}
 x_t= cdt + {\sqrt \alpha}dB_t,
\end{equation}
where  $\alpha$ is volatility.
This equation can be integrated to be
\begin {equation}
 x_t=\int _0^t{c(t-s)}\sqrt\alpha dB_s .
\end{equation}
Let
\begin {equation}
 X_t=X_0 e^{x_t},.
\end{equation}
and applying Ito formula we obtain
\begin {equation}
d X_t=\mu X_t dt + {\sqrt \alpha}X_tdB_t,
\end{equation}
where
let
\begin {equation}
\mu = c+ \alpha,
\end{equation}
which represent averaged income,
and this can be integrated to be
\begin {equation}
 X_t=e^{{\sqrt \alpha} B_t+{{(\mu - \alpha)\over 2}{t }}}.
\end{equation}
This is the optimal pricing of European call option known as Black-Scholes formula.$^1$ It will be worth to note that eq.(1.6) include perturbed drift as a result there is an underlying displacement.

\subsection{Trading Strategy}
\par
Trading strategy is a twin of stochastic variables denoted by$^2$
\begin {equation}
 (a_t,b_t).
\end{equation}
In this equation $a_t$ and $b_t$ are given by as defined in the space of square integral function $L^2$
\begin {equation}
a_tdX_t=\mu X_tdt +{ \sqrt\alpha}X_tdB_t
\end{equation}
and
\begin {equation}
db_t=rb_t\beta_t dt, 
\end{equation}
where $r$ denotes interest ratio and $\beta_t$ is pricing process of bond given by
\begin {equation}
d\beta_t=rb_t\beta_t dt.
\end{equation}
Using these notation the portfolio is given by
 \begin {equation}
V_t=a_tX_t+b_t\beta_t.
\end{equation}

\subsection{Adaptive Optimization}
\par
Let us consider a stochastic integral of optimal pricing  w.r.t.  trading strategy as follows
\begin {equation}
 X_t=\int^t_0 X_s d\{a_s,b_s\}_s .
\end{equation}
If r.h.s. does not change pricing is called adoptive to trading strategy. This defines a  martingale.$^3$
This means that pricing is globally optimized to trading strategy. This equation can be solved by iteration as follows.
\par 
{\hskip 1 in \bf repeat until}
\begin {equation}
{\bf E(X_{m}|\{a_n,b_n\}_n)= X_n < \epsilon,  \, a.e. \, \,  m>n}
\end{equation}
\par
{\hskip 1 in \bf continue}
\par
This provides an algorithm implementable into digital computer where  the optimal pricing adoptive to trading strategy in which r.h.s denotes a conditional probabilty. In this case, trading strategy needs not be self-financing.

\section{Concluding remarks}
\par
Optimal pricing for risky asset adoptive to trading strategy was formulated in terms of Ito type stochastic differential equation. It should be noted, however, that there is a limitation in present formulation that it rely on Brownian motion. As a result it depends on Gaussian process based on central limit theorem which states summation of independent variable is governed by Gaussian distribution. The shortcomings arises from the fact that variables are not independent in the real world.



\bibliographystyle{amsplain}

\begin{thebibliography}{10}

\bibitem {1} T.Mikosh, \textit{Elementary Stochastic Calculus with Finance in View}  (World Scientific Pub. 1998).

\bibitem {2} D.Duffie, \textit{Dynamic Asset Pricing Theory} (Princeton University Press 1996).

\bibitem {3} T.Hida, \textit{Brownian Motion} (Springer-Verlag 1980).

\end{thebibliography}

\end{document}